\documentclass[reqno,psamsfonts]{amsart}
\usepackage{amssymb}
\usepackage{lattice}
\usepackage{amscd}
\usepackage{upref}
\usepackage{graphics}

\theoremstyle{plain}
\newtheorem*{theoremn}{Theorem}
\newtheorem{theorem}{Theorem}
\newtheorem{lemma}{Lemma}[section]
\newtheorem{proposition}[lemma]{Proposition}

\newtheorem*{stat}{\name}
\newcommand{\name}{testing}
\newcommand{\Pow}[1]{\operatorname{P}(#1)} %power set

\theoremstyle{definition}
\newtheorem*{definition}{Definition}

\newtheorem{problem}{Problem}

\newcommand{\ootimes}{\mathbin{\ol{\otimes}}}
\newcommand{\ootimest}{\mathbin{\vec{\otimes}}}
\newcommand{\jz}{$\set{\jj,0}$}
\newcommand{\tS}{\widetilde{S}}
\newcommand{\id}{\mathrm{id}}

\DeclareMathOperator{\J}{J}
\DeclareMathOperator{\Hom}{Hom}
\DeclareMathOperator{\Id}{Id}

\begin{document}

\title{Flat semilattices}

\author{G.~Gr\"atzer}
\thanks{The research of the first author was partially
        supported by the NSERC of Canada.}
\address{Department of Mathematics\\
         University of Manitoba\\
         Winnipeg, Manitoba\\
         Canada R3T 2N2}
\email{gratzer@cc.umanitoba.ca}
 \urladdr{http://www.maths.umanitoba.ca/homepages/gratzer.html/}

\author{F.~Wehrung}
\address{C.N.R.S.\\
         D\'epartement de Math\'ematiques\\
         Universit\'e de Caen\\
         14032 Caen Cedex\\
         France}
\email{gremlin@math.unicaen.fr}
\urladdr{http://www.math.unicaen.fr/\~{}wehrung}

\date{June 1, 1998}

\begin{abstract}
Let $A$, $B$, and $S$ be \jz-semilattices and let $f \colon
A \hookrightarrow B$ be a \jz-semilattice embedding.  Then the canonical
map, $f \otimes \id_S$, of the tensor product $A \otimes S$ into the
tensor product $B \otimes S$ is not necessarily an embedding.

The \jz-semilattice $S$ is \emph{flat}, if for every embedding
$f \colon A \hookrightarrow B$, the canonical map $f\otimes\id$ is
an embedding. We prove that \emph{a \jz-semi\-lat\-tice $S$ is flat if
and only if it is distributive}. 
\end{abstract}

\maketitle
\section*{Introduction}
Let $A$ and $B$ be \jz-semilattices. We denote by $A \otimes B$ the
\emph{tensor product} of $A$ and $B$, defined as the free \jz-semilattice
generated by the set 
 \[
   (A - \set{0}) \times (B - \set{0})
 \]
 subject to the relations 
 \[
   \vv<a,b_0> \jj \vv<a,b_1> = \vv<a, b_0 \jj b_1>,
 \]
 for $a  \in A - \set{0}$, $b_0$, $b_1  \in B - \set{0}$; and
symmetrically, 
 \[
   \vv<a_0, b> \jj \vv<a_1, b> = \vv<a_0 \jj a_1, b>,
 \]
 for $a_0$, $a_1  \in  A - \set{0}$,  $b  \in B - \set{0}$.

$A \otimes B$ is a universal object with respect to a natural notion of
\emph{bimorphism}, see \cite{Fras78}, \cite{GLQu81}, and \cite{GrWe}.
This definition is similar to the classical definition of the tensor
product of modules over a commutative ring. Thus, for instance,
\emph{flatness} is defined similarly: The \jz-semilattice $S$ is
\emph{flat}, if for every embedding $f \colon A \hookrightarrow B$, the
canonical map $f \otimes \id_S \colon A \otimes S \to B \otimes S$ is an
embedding.

Our main result is the following: 

\begin{theoremn}
Let $S$ be a \jz-semilattice. Then $S$ is flat if{f} $S$ is distributive.
\end{theoremn}

\section{Background}\label{S:antitone}

\subsection{Basic concepts}\label{S:basic}
We shall adopt the notation and terminology of \cite{GrWe}. In
particular, for every \jz-semilattice $A$, we use the notation $A^-=A -
\set{0}$.  Note that $A^-$ is a subsemilattice of $A$.

A semilattice $S$ is \emph{distributive}, if whenever $a \leq b_0 \jj
b_1$ in $S$, then there exist $a_0 \leq b_0$ and $a_1 \leq b_1$ such
that $a = a_0 \jj a_1$; equivalently, if{f} the lattice $\Id S$ of all
ideals of $S$, ordered under inclusion, is a distributive lattice; see
\cite{Grat78}.

\subsection{The set representation}\label{S:set}
In \cite{GrWe}, we used the following representation of the tensor product.

First, we introduce the notation:
 \[
   \nabla_{A,B} = (A \times \set{0})  \uu  (\set{0} \times B).
 \]

Second, we introduce a partial binary operation on $A \times B$: let
$\vv<a_0,b_0>$, $\vv<a_1,b_1> \in A \times B$; the \emph{lateral join} of
$\vv<a_0,b_0>$ and $\vv<a_1,b_1>$ is defined if $a_0 = a_1$ or $b_0 = b_1$,
in which case, it is the join, $\vv<a_0\jj a_1,b_0\jj b_1>$.

Third, we define bi-ideals: a nonempty subset $I$ of $A \times B$ is a
\emph{bi-ideal} of $A \times  B$, if it satisfies the following conditions:
 \begin{enumerate}
 \item $I$ is hereditary;
 \item $I$ contains $\nabla_{A,B}$;
 \item $I$ is closed under lateral joins.
 \end{enumerate}

The \emph{extended tensor product} of $A$ and $B$, denoted by $A \ootimes
B$, is the lattice of all bi-ideals of $A \times B$.

It is easy to see that $A \ootimes  B$ is an algebraic lattice. For
$a \in A$ and $b \in B$, we define $a \otimes b \in A \ootimes B$ by
 \[
   a \otimes b =\nabla_{A,B}  \uu  \setm{\vv<x, y>  \in A \times B}
  {\vv<x, y>  \leq \vv<a, b>}
 \]
and call $a \otimes b$ a \emph{pure tensor}.  A pure tensor is a
principal (that is, one-generated) bi-ideal.  

Now we can state the representation:

\begin{proposition}\label{P:TensSet}
The tensor product $A \otimes B$
can be represented as the
\jz-sub\-semi\-lat\-tice of compact elements of $A \ootimes B$.
\end{proposition}

\subsection{The construction of $A\ootimest B$}\label{S:hom}
 The proof of the Theorem uses the following representation of the
tensor product, see J. Anderson and N. Kimura \cite{AnKi78}.

 Let $A$ and $B$ be \jz-semilattices. Define
 \[
   A \ootimest B = \Hom(\vv<A^-; \jj>,\vv<\Id B; \ii>)
 \]
 and for $\gx \in A \ootimest B$, let
 \[
  \ge(\gx) = \setm{\vv<a, b> \in A^{-} \times B^{-}}
     {b \in \gx(a)} \uu \nabla_{A,B}.
 \]

\begin{proposition}\label{P:TensAnti}
 The map $\ge$ is an order preserving isomorphism between $A \ootimes B$
and $A\ootimest B$ and, for $H \in A \ootimes B$, $\ge^{-1}(H)$ is given
by the formula
 \[
   \ge^{-1}(H)(a) = \setm{b \in B}{\vv<a, b> \in H},
 \]
 for $a \in A^-$.
\end{proposition}

If $a\in A$ and $b\in B$, then $\ge(a\otimes b)$ is the map
$\gx\colon A^-\to\Id B$:
\[
\gx(x)=
   \begin{cases}
    (b], &\text{if $x\leq a$};\\
   \set{0}, &\text{otherwise.}
   \end{cases}
\]

If $A$ is \emph{finite}, then a homomorphism from $\vv<A^-;\jj>$ to
$\vv<\Id B;\ii>$ is determined by its restriction to $\J(A)$, the set of
all join-irreducible elements of $A$.  For example, let $A$ be a finite
Boolean semilattice, say $A= \Pow n$ ($n$ is a non-negative integer,
$n=\set{0, 1, \ldots, n - 1}$), then $A \ootimes B \iso (\Id B)^n$, and
the isomorphism from $A \ootimes B$ onto $(\Id B)^n$ given by
Proposition~\ref{P:TensAnti} is the unique complete \jz-homomorphism
sending every element of the form $\set{i} \otimes b$ ($i < n$ and $b
\in B$) to $\vv<(\gd_{ij}b] \mid j < n>$ (where $\gd_{ij}$ is the
Kronecker symbol). If $n = 3$, let $\gb \colon \Pow 3 \ootimes
S \to (\Id S)^3$ denote the natural isomorphism.

Next we compute $A \ootimest B$, for $A= {M}_3$, the diamond, and
$A={N}_5$, the pentagon (see Figure 1). In the following two subsections,
let $S$ be a \jz-semilattice. Furthermore, we shall denote by $\tS$ the
ideal lattice of $S$, and identify every element $s$ of $S$ with its
image, $(s]$, in $\tS$. 

\subsection{The lattices ${M}_3\ootimes S$ and ${M}_3[\tS]$; the map
$i$.}

\begin{figure}[hbt]
\centerline{\includegraphics{IIs1f1.ill}}
\end{figure}

Let ${M}_3 = \set{0, p, q, r, 1}$, $\J({M}_3) = \set{p, q, r}$ (see
Figure~1). The nontrivial relations of $\J({M}_3)$ are the
following:
 \begin{equation}\label{Eq:CovM3}
   p < q \jj r,\quad q < p \jj r,\quad\text{and}\quad r < p \jj q.
 \end{equation} Accordingly, for every lattice $L$, we define
 \begin{equation}\label{Eq:RepM3}
   {M}_3[L] = \setm{\vv<x, y, z>\in L^3}{x \mm y = x \mm z = y \mm z}
 \end{equation} (this is the \emph{Schmidt's construction}, see
\cite{Quac85} and \cite{Schm68}). The isomorphism from ${M}_3\ootimes S$
onto ${M}_3[\tS]$ given by Proposition~\ref{P:TensAnti} is the unique
complete \jz-homomorphism $\ga$ such that, for all $x\in S$, 
 \begin{align*}
   \ga(p\otimes x) &= \vv<x,0,0>,\\
   \ga(q\otimes x) &= \vv<0,x,0>,\\
   \ga(r\otimes x) &= \vv<0,0,x>.
 \end{align*}

We shall make use later of the unique \jz-embedding
 \[
   i\colon {M}_3 \hookrightarrow \Pow 3
 \]
 defined by
 \begin{align*}
   i(p) &=\set{1,2},\\
   i(q) &=\set{0,2},\\
   i(r) &=\set{0,1}.
 \end{align*}

\subsection{The lattices ${N}_5\ootimes S$ and
${N}_5[\tS]$; the map $i'$.}

Let ${N}_5 = \set{0, a, b, c, 1}$, $\J({N}_5) = \set{a, b, c}$ with
$a > c$ (see Figure~1). The nontrivial relations of $\J({N}_5)$ are the
following:
 \begin{equation}\label{Eq:CovN5}
   c < a \quad \text{and} \quad a < b \jj c.
 \end{equation}
 Accordingly, for every lattice $L$, we define
 \begin{equation}\label{Eq:RepN5}
   {N}_5[L] = \setm{\vv<x, y, z> \in L^3}{y \mm z \leq x \leq z}.
 \end{equation}
 The isomorphism from ${N}_5\ootimes S$ onto ${N}_5[\tS]$, given by
Proposition~\ref{P:TensAnti}, is the unique complete
\jz-homomorphism $\ga'$ such that, for all $x \in S$, 
 \begin{align*}
   \ga'(a\otimes x) &=\vv<x,0,x>,\\
   \ga'(b\otimes x) &=\vv<0,x,0>,\\
   \ga'(c\otimes x) &=\vv<0,0,x>.
 \end{align*}
 We shall make use later of the unique \jz-embedding
 \[
   i'\colon {N}_5 \hookrightarrow \Pow 3
 \]
 defined by
  \begin{align*}
   i'(a) &= \set{0, 2},\\
   i'(b) &= \set{1, 2},\\
   i'(c) &= \set{0}.
 \end{align*}
 
\subsection{The complete homomorphisms $f \ootimes g$}
The proof of the following lemma is straightforward:

 \begin{lemma}\label{L:ooemb}
 Let $A$, $B$, $A'$, and $B'$ be \jz-semilattices, let $f\colon A\to A'$
and $g\colon B\to B'$ be \jz-homomorphisms. Then the natural
\jz-homomorphism $h = f \otimes g$ from $A \otimes B$ to $A' \otimes B'$
extends to a unique complete \jz-homomorphism $\ol{h} = f\ootimes g$ from
$A \ootimes B$ to $A' \ootimes B'$. Furthermore, if $h$ is an embedding,
then $\ol{h}$ is also an embedding.
 \end{lemma}

We refer to Proposition 3.4 of \cite{GrWe} for an explicit description
of the map $\ol{h}$.

\section{Characterization of flat \jz-semilattices}

Our definition of flatness is similar to the usual one for modules over
a commutative ring:

 \begin{definition}
 A \jz-semilattice $S$ is \emph{flat}, if for every embedding $f \colon
A \hookrightarrow B$ of \jz-semilattices, the tensor map
$f \otimes \id_S \colon A \otimes S \to B \otimes S$ is an embedding.
 \end{definition}

In this definition, $\id_S$ is the identity map on $S$.

In Lemmas~\ref{L:CD}--\ref{L:NoN5}, let $S$ be a \jz-semilattice and we
assume that both homomorphisms $f = i \otimes \id_S$ and
$f' = i' \otimes \id_S$ are embeddings.

As in the previous section, we use the notation $\tS = \Id S$, and
identify every element $s$ of $S$ with the corresponding principal ideal
$(s]$.

We define the maps $g \colon {M}_3[\tS] \to \tS^3$ and
$g' \colon {N}_5[\tS] \to \tS^3$ by the following formulas:
 \begin{align*}
\text{For all }&\vv<x,y,z>\in {M}_3[\tS],&
g(\vv<x,y,z>)&=\vv<y\jj z,x\jj z,x\jj y>,\\
\text{For all }&\vv<x,y,z>\in {N}_5[\tS],&
g'(\vv<x,y,z>)&=\vv<z,y,x\jj y>.
 \end{align*}

Note that $g$ and $g'$ are complete \jz-homomorphisms. The proof of the
following lemma is a straightforward calculation.

\begin{lemma}\label{L:CD}
The following two diagrams commute:
\[
 \begin{CD}
 {M}_3\ootimes S @>f>> \Pow 3\ootimes S\\
 @V{\ga}VV @VV{\gb}V\\
 {M}_3[\tS] @>>g> \tS^3
 \end{CD}
\]
\[
\begin{CD}
 {N}_5\ootimes S @>f'>> \Pow 3\ootimes S\\
 @V{\ga'}VV @VV{\gb}V\\
 {N}_5[\tS] @>>g'> \tS^3
 \end{CD}
\]
Therefore, both $g$ and $g'$ are embeddings.
\end{lemma}

\begin{lemma}\label{L:NoM3}
The lattice $\tS$ does not contain a copy of ${M}_3$.
\end{lemma}

\begin{proof}
Suppose, on the contrary, that $\tS$ contains a copy of
${M}_3$, say $\set{o, x, y, z, i}$ with $o < x$, $y$, $z < i$. Then
both elements $u=\vv<x, y, z>$ and $v=\vv<i, i, i>$ of $L^3$ belong to
${M}_3[\tS]$, and $g(u) = g(v) = \vv<i, i, i>$. This contradicts the
fact, proved in Lemma~\ref{L:CD}, that $g$ is one-to-one.
\end{proof}

 \begin{lemma}\label{L:NoN5}
The lattice $\tS$ does not contain a copy of ${N}_5$.
 \end{lemma}

\begin{proof}
Suppose, on the contrary, that $\tS$ contains a copy of
${N}_5$, say $\set{o, x, y, z, i}$ with $o < x < z <i$ and $o < y < i$.
Then both elements $u = \vv<x, y, z>$ and $v=\vv<z ,y, z>$ of $L^3$
belong to ${N}_5[\tS]$, and $g'(u) = g'(v) = \vv<z ,y, i>$. This
contradicts the fact, proved in Lemma~\ref{L:CD}, that $g'$ is
one-to-one.
\end{proof}

Lemmas \ref{L:NoM3} and \ref{L:NoN5} together prove that $\tS$ is
distributive, and therefore $S$ is a distributive semilattice. Now we
are in position to prove the main result of this paper in the following
form:

\begin{theorem}\label{T:FlatSem}
 Let $S$ be a \jz-semilattice. Then the following are
equivalent:
 \begin{enumerate}
 \item $S$ is flat.
 \item Both homomorphisms $i\otimes\id_S$ and $i'\otimes\id_S$ are
embeddings.
 \item $S$ is distributive.
 \end{enumerate}
 \end{theorem}

\begin{proof}\hfill

(i) \emph{implies} (ii).  This is trivial.

(ii) \emph{implies} (iii). This was proved in Lemmas \ref{L:NoM3}
and \ref{L:NoN5}.

 (iii) \emph{implies} (i). Let $S$ be a distributive \jz-semilattice; we
prove that $S$ is flat. Since the tensor product by a fixed factor
preserves direct limits (see Proposition~2.6 of \cite{GrWe}), flatness
is preserved under direct limits. By P.~Pudl\'ak~\cite{Pudl85}, every
distributive join-semilattice is the direct union of all its finite
distributive subsemilattices; therefore, it suffices to prove that every
finite distributive \jz-semilattice $S$ is flat. Since $S$ is a
distributive lattice, it admits a lattice embedding into a finite
Boolean lattice $B$. We have seen in Section~\ref{S:hom} that if
$B = \Pow n$, then $A\otimes B = A^n$ (up to a natural isomorphism),
for every \jz-semilattice $A$.  It follows that $B$ is~flat. Furthermore,
the inclusion map $S\hookrightarrow B$ is a lattice embedding; in
particular, with the terminology of \cite{GrWe}, an
\emph{L-homomorphism}. Thus, the natural map from $A\otimes S$ to
$A\otimes B$ is, by Proposition 3.4 of \cite{GrWe}, a \jz-semilattice
embedding. This implies the flatness of $S$.
 \end{proof}

\section{Discussion}
 It is well-known that a module over a given principal ideal
domain $R$ is flat if and only if it is torsion-free, which is
equivalent to the module being a direct limit of (finitely generated)
free modules over $R$. So the analogue of the concept of torsion-free
module for semilattices is be the concept of distributive
semilattice. This analogy can be pushed further, by using the following
result, proved in \cite{GoWe}: \emph{a~join-semilattice is distributive
if{}f it is a direct limit of finite Boolean semilattices}.

 \begin{problem}
 Let $\mathbf{V}$ be a variety of lattices. Let us say that a
\jz-semilattice $S$ is \emph{in $\mathbf{V}$}, if $\Id S$ as a lattice
is in $\mathbf{V}$. Is every \jz-semilattice in $\mathbf{V}$ a direct
limit (resp., direct union) of \emph{finite} join-semilattices in
$\mathbf{V}$?
 \end{problem}

If $\mathbf{V}$ is the variety of all lattices, we obtain the obvious
result that every \jz-semilattice is the direct union of its finite
\jz-subsemilattices. If $\mathbf{V}$ is the variety of all distributive
lattices, there are two results (both quoted above): P. Pudl\'ak's result
and K. R. Goodearl and the second author's result.

 \begin{problem}
 Let $\mathbf{V}$ be a variety of lattices. When is a \jz-semilattice
$S$ flat with respect to \jz-semilattice embeddings in $\mathbf{V}$?
That is, when is it the case that for all \jz-semilattices $A$ and $B$
in $\mathbf{V}$ and every semilattice embedding $f\colon
A\hookrightarrow B$, the natural map $f\otimes\id_S$ is an embedding?
 \end{problem}

\end{document}